\documentclass[12pt]{article}

\usepackage[margin=1.3in, top=1.3in, bottom=1.3in]{geometry}

\usepackage[utf8]{inputenc} 
\usepackage[T1]{fontenc}    
\usepackage{url}            
\usepackage{booktabs}       
\usepackage{amsfonts}       
\usepackage{nicefrac}       
\usepackage{microtype}      

\usepackage{amsmath,amssymb,amsthm,xcolor}
\usepackage{array}
\usepackage{float}

\newtheorem{theorem}{Theorem}
\newtheorem{lemma}{Lemma}

\newtheorem{definition}{Definition}
\usepackage{graphicx}
\usepackage{subcaption}
\usepackage{mdframed}
\usepackage{lipsum}
\usepackage{wrapfig}
\usepackage[hidelinks,hypertexnames=false]{hyperref}
\usepackage{tikz}
\usepackage{enumitem}
\allowdisplaybreaks
\usepackage{pgfplots}
\usetikzlibrary{intersections, pgfplots.fillbetween}
\usepackage{epstopdf}
\usepackage{algorithmic}
\usepackage{booktabs}
\usepackage{changepage}

\title{\LARGE Lyapunov stability of the subgradient method with constant step size}

\begin{document}

\author{\large C\'edric Josz\thanks{\url{cj2638@columbia.edu}, IEOR, Columbia University, New York. Research supported by NSF EPCN grant 2023032 and ONR grant N00014-21-1-2282.} \and Lexiao Lai\thanks{\url{ll3352@columbia.edu}, IEOR, Columbia University, New York.}}
\date{}

\maketitle
\begin{center}
    \textbf{Abstract}
    \end{center}
    \vspace*{-3mm}
 \begin{adjustwidth}{0.2in}{0.2in}
~~~~We consider the subgradient method with constant step size for minimizing locally Lipschitz semi-algebraic functions. In order to analyze the behavior of its iterates in the vicinity of a local minimum, we introduce a notion of discrete Lyapunov stability and propose necessary and sufficient conditions for stability.
\end{adjustwidth} 
\vspace*{3mm}


\section{Introduction}
\label{sec:Introduction}
The subgradient method with constant step size for minimizing a locally Lipschitz function $f:\mathbb{R}^n\rightarrow\mathbb{R}$ consists in choosing an initial point $x_0 \in \mathbb{R}^n$ and generating a sequence of iterates according to the update rule $x_{k+1} \in x_k - \alpha \partial f(x_k), ~ \forall k \in \mathbb{N} := \{0,1,2,\hdots\}$, where $\alpha>0$ is the step size and $\partial f$ is the Clarke subdifferential \cite[Chapter 2]{clarke1990}. While this method is often used in practice to solve nonconvex and nonsmooth problems, there is little theoretical understanding of the behavior of its iterates. To the best of our knowledge, the only known results are in the convex setting, as we next describe.

If $f$ is convex and the Euclidean norm of its subgradients is bounded above by a constant $c$, then $\liminf f(x_k) -\inf f \leqslant c^2\alpha/2$ provided that the infimum is reached \cite[Proposition 3.2.3]{bertsekas2015convex}. In order to get within $\epsilon$ accuracy of that bound, $\lfloor d(x_0,X)^2/(\alpha \epsilon) \rfloor$ iterations suffice where $\lfloor \cdot \rfloor$ denotes the floor of a real number and $d(x_0,X)$ is the distance between the initial iterate $x_0$ and the set of minimizers $X \subset \mathbb{R}^n$ \cite[Proposition 3.2.4]{bertsekas2015convex}. If the objective function grows quadratically (at least as fast as $t\in \mathbb{R} \mapsto\beta t^2$ for some $\beta>0$) around the set of minimizers, then the iterates asymptotically get within $c\sqrt{\alpha}/\sqrt{2\beta}$ distance to the set of minimizers if $\alpha \in (0,1/(2\beta)]$ \cite[Proposition 3.2.5]{bertsekas2015convex}. If we relax the boundedness assumption on the subgradients to $\|s\| \leqslant c \sqrt{1+d(x,X)^2}$ for all $(x,s)$ in the graph of $\partial f$, then we get the slightly weaker bound $\liminf f(x_k) -\inf f \leqslant c^2\alpha/2(1+d(x_0,X))$ \cite[Exercise 3.6]{bertsekas2015convex}.

Given the absence of theoretical results in the nonconvex setting, in this note we take a first step by investigating the behavior of the subgradient method in the vicinity of a local minimum of $f$. In order to do so, we propose a notion of stability akin to Lyapunov stability in dynamical systems \cite{liapounoff1907probleme} \cite[Equation (5.6)]{sastry2013nonlinear}. Informally, a point is stable if all of the iterates of the subgradient method remain in any neighborhood of it, provided that the initial point is close enough to it and that the step size is small enough.

Without any further assumptions on $f$, the notions of stability and local optimality are decorrelated. Indeed, the classical counterexample $f(x) = x^2\sin(1/x)$ admits a stable point which is not a local minimum, while the Rockafellar function \cite[p. 5]{rockafellar1981favorable} (see also \cite[Proposition (1.9)]{lebourg1979generic}) admits a local minimum that is not stable. Assuming sharpness \cite[Assumption A 2]{davis2018subgradient} and weak convexity \cite[Assumption A 1]{davis2018subgradient}, it can easily be shown that strict local minima are stable. These assumptions may not hold in practice however \cite[p. 121]{davis2020stochastic} \cite[2.3.6 Proposition]{clarke1990}, and can be difficult to check \cite[Conjecture 8.7]{charisopoulos2021low}. We thus confine our investigation to locally Lipschitz semi-algebraic functions \cite{bochnak2013real,pham2016genericity}. The tame generalization \cite{van1996geometric} is immediate and captures seemingly all applications of interest nowadays.

When $f$ is locally Lipschitz and semi-algebraic, the set of stable points and local minima coincide in two cases. The first is when $f$ is continuously differentiable with a locally Lipschitz gradient. The fact that local minima are stable in this regime can be deduced using arguments from \cite[Proposition 3.3]{absil2005convergence}. The converse is a consequence of one of our results (Theorem \ref{thm:necessary}). The second case is that of continuous-time subgradient dynamics, where instead of iterates one considers absolutely continuous solutions to the differential inclusion $x' \in -\partial f(x)$. This is a simple generalization of \cite[Theorem 3]{absil2006stable} which holds for real analytic functions. Much of modern numerical optimization however falls outside the scope of these two cases, namely that of smooth objective functions and continuous-time dynamics. It is thus important to determine in the discrete case and when $f$ is merely locally Lipschitz and semi-algebraic, whether local minima are stable, and conversely, whether stable points are local minima. In this note, we show that for a point to be stable, it is necessary for it to be a local minimum and it suffices for it to be a strict local minimum. 

\section{Characterizing stability}
\label{sec:Characterizing stability}

Let $\|\cdot\|$ be the induced norm of an inner product $\langle \cdot, \cdot\rangle$ on $\mathbb{R}^n$. Let $B(a,r)$ and $\mathring{B}(a,r)$ respectively denote the closed ball and the open ball of center $a\in \mathbb{R}^n$ and radius $r>0$. We next define the notion of discrete Lyupanov stability.

\begin{definition}
\label{def:discrete_lyapunov}
\normalfont
We say that $x^*\in \mathbb{R}^n$ is a stable point of a locally Lipschitz function $f:\mathbb{R}^n\rightarrow\mathbb{R}$ if for all $\epsilon>0$, there exist $\delta>0$ and $\bar{\alpha}>0$ such that for all $\alpha \in (0,\bar{\alpha}]$, the subgradient method with constant step size $\alpha$ initialized in $B(x^*,\delta)$ has all its iterates in $B(x^*,\epsilon)$.
\end{definition}

In order to characterize stability, we rely on the theory of differential inclusions \cite{aubin1984differential}. However, existing results cannot directly be applied, mainly because the locally Lipschitz assumption is too weak. We thus adapt them to our setting via Lemma \ref{lemma:euler_method}. Note that the analysis of the stochastic subgradient method also involves differential inclusions \cite{benaim2005stochastic,bianchi2019constant,davis2020stochastic,bolte2020conservative,bianchi2022convergence}.

\begin{lemma}
\label{lemma:euler_method}
    Let $F:\mathbb{R}^n\rightrightarrows\mathbb{R}^n$ be an upper semicontinuous mapping with nonempty, compact, and convex values. Let $X_0$ be a compact subset of $\mathbb{R}^n$ and let $T>0$. Assume that there exist $\bar{\alpha},r>0$ such that for all $\alpha \in (0,\bar{\alpha}]$ and for all sequence $(x_k)_{k\in\mathbb{N}}$ such that
    \begin{equation}
    \label{eq:Euler}
                x_{k+1} \in x_k + \alpha F(x_k), ~~~ \forall k \in \mathbb{N},~~~ x_0 \in X_0,
    \end{equation}
    we have that $x_0,\hdots,x_{\lfloor T/\alpha \rfloor+1} \in B(0,r)$. For all $\epsilon>0$, there exists $\hat{\alpha}\in (0,\bar{\alpha}]$ such that for all $\alpha \in (0,\hat{\alpha}]$ and for all sequence $(x_k)_{k\in\mathbb{N}}$ satisfying \eqref{eq:Euler}, there exists an absolutely continuous function $x:[0,T]\rightarrow\mathbb{R}^n$ such that
    \begin{equation}
    \label{eq:ivp}
                x'(t) \in F(x(t)),~~~ \text{for a.e.}~ 
                t \in (0,T), ~~~ x(0) \in X_0, 
    \end{equation}
    for which $\|\bar{x}(t) - x(t) \| \leqslant \epsilon$ for all $t \in [0,T]$, where $\bar{x}:[0,T]\rightarrow\mathbb{R}^n$ is defined by $\bar{x}(t) := x_k + (t-\alpha k)/\alpha  (x_{k+1}-x_k)$ for all $t\in [\alpha k , \min\{\alpha(k+1),T\}]$ and $k\in \{0,\hdots,\lfloor T/\alpha\rfloor\}$.
\end{lemma}

\begin{proof}
Let $(\alpha_m)_{m \in \mathbb{N}}$ denote a sequence of positive numbers that converges to zero. Without loss of generality, we may assume that the sequence is bounded above by $\bar{\alpha}$. To each term in the sequence, we attribute a sequence $(x_k^{m})_{k\in\mathbb{N}}$ generated by the Euler method with step size $\alpha_m$ and initialized in $X_0$, that is to say, which satisfies \eqref{eq:Euler} with $\alpha:=\alpha_m$. By assumption, $x_0^{m},\hdots,x_{\lfloor T/\alpha_m \rfloor+1}^{m} \in B(0,r)$. Consider the linear interpolation of those iterates, that is to say, the function $\bar{x}^{m}:[0,T]\rightarrow\mathbb{R}^n$ defined by $\bar{x}^{m}(t) := x_k^{m} + (t-\alpha_m k)(x_{k+1}^{m}-x_k^{m})/\alpha_m$ for all $t\in [\alpha_m k , \min\{\alpha_m(k+1),T\}]$ and $k\in \{0,\hdots,\lfloor T/\alpha_m \rfloor\}$. Since $B(0,r)$ is convex, it holds that $\|\bar{x}^{m}(t)\| \leqslant r $ for all $t\in [0,T]$. In addition, since $F$ is upper semicontinuous and compact valued, by \cite[Proposition 3 p. 42]{aubin1984differential} there exists $r'>0$ such that $F(B(0,r)) \subset B(0,r')$. Observe that $(\bar{x}^{m})'(t) = (x_{k+1}^{m}-x_k^{m})/\alpha_m \in F(x_k^{m})$ for all $t\in (\alpha_m k , \min\{\alpha_m(k+1),T\})$ and $k\in \{0,\hdots,\lfloor T/\alpha_m \rfloor\}$. Hence, we have that $\|(\bar{x}^{m})'(t)\| \leqslant r'$ for almost every $t\in (0,T)$. By successively applying the Arzel\`a-Ascoli and the Banach-Alaoglu theorems (see \cite[Theorem 4 p. 13]{aubin1984differential}), there exists a subsequence (again denoted $(\alpha_{m})_{m\in \mathbb{N}}$) and an absolutely continuous function $x:[0,T]\rightarrow \mathbb{R}^n$ such that $\bar{x}^{m}(\cdot)$ converges uniformly to $x(\cdot)$ and $(\bar{x}^{m})'(\cdot)$ converges weakly to $x'(\cdot)$ in $L^1([0,T],\mathbb{R}^n)$. Furthermore, for all $t\in (\alpha_{m} k , \min\{\alpha_{m}(k+1),T\})$ and $k \in \{0,\hdots,\lfloor T/\alpha_m \rfloor\}$, observe that
\begin{subequations}
\begin{align*}
    (\bar{x}^{m}(t),(\bar{x}^{m})'(t)) 
    & = \left(x_k^{m} + (t-\alpha_m k)\frac{x_{k+1}^{m}-x_k^{m}}{\alpha_m},\frac{x_{k+1}^{m}-x_k^{m}}{\alpha_m}\right) \\
    & \in \left(\{x_k^{m}\} + (t-\alpha_m k)F(x_k^{m})\right) \times F(x_k^{m})  \\
    & = \{x_k^{m}\} \times F(x_k^{m}) + (t-\alpha_m k)F(x_k^{m}) \times \{0\}  \\
    & \subset \mathrm{graph}(F) + B(0,r'\alpha_m)\times \{0\}.
\end{align*}
\end{subequations}
According to \cite[Convergence Theorem p. 60]{aubin1984differential}, it follows that $x'(t) \in F(x(t))$ for almost every $t\in (0,T)$. The sequence of initial points $(x^{m}(0))_{m\in\mathbb{N}}$ lies in the closed set $X_0$, hence its limit $x(0)$ lies in $X_0$ as well. As a result, $x(\cdot)$ is a solution to the differential inclusion \eqref{eq:ivp}. 

To sum up, we have shown that for every sequence $(\alpha_m)_{m \in \mathbb{N}}$ of positive numbers converging to zero, there exists a subsequence for which the corresponding linear interpolations uniformly converge towards a solution of the differential inclusion \eqref{eq:ivp}. The conclusion of the theorem now easily follows. To see why, one can reason by contradiction and assume that there exists $\epsilon>0$ such that for all $\hat{\alpha} \in (0,\bar{\alpha}]$, there exist $\alpha \in (0,\hat{\alpha}]$ and a sequence $(x_k)_{k \in \mathbb{N}}$ generated by the Euler method with step size $\alpha$ and initialized in $X_0$ such that, for any solution $x(\cdot)$ to the differential inclusion \eqref{eq:ivp}, it holds that $\|\bar{x}(t) - x(t) \| > \epsilon$ for some $t \in [0,T]$. We can then generate a sequence $(\alpha_m)_{m \in \mathbb{N}}$ of positive numbers converging to zero such that, for any solution $x(\cdot)$ to the differential inclusion \eqref{eq:ivp}, it holds that $\|\bar{x}^{m}(t) - x(t) \| > \epsilon$ for some $t \in [0,T]$. Since there exists a subsequence $(\alpha_{\varphi(m)})_{m \in \mathbb{N}}$ such that $\bar{x}^{\varphi(m)}$ uniformly converges to a solution to the differential inclusion \eqref{eq:ivp}, we obtain a contradiction. 
\end{proof}

We will apply Lemma \ref{lemma:euler_method} to the case where the set-valued mapping $F$ is the opposite of the Clarke subdifferential of a locally Lipschitz function $f:\mathbb{R}^n\rightarrow\mathbb{R}$. Recall that a point $x^* \in \mathbb{R}^n$ is a local minimum (respectively strict local minimum) of a function $f:\mathbb{R}^n\rightarrow\mathbb{R}$ if there exists a positive constant $\epsilon$ such that $f(x^*) \leqslant f(x)$ for all $x \in B(x^*,\epsilon)\setminus \{x^*\}$ (respectively $f(x^*) < f(x)$). Using Lemma \ref{lemma:euler_method}, we obtain the following necessary condition for stability.

\begin{theorem}
\label{thm:necessary}
Stable points of locally Lipschitz semi-algebraic functions are local minima.
\end{theorem}

\begin{proof}
Let $x^* \in \mathbb{R}^n$ denote a stable point of a locally Lipschitz semi-algebraic function $f:\mathbb{R}^n\rightarrow\mathbb{R}$. We reason by contradiction and assume that $x^*$ is not a local minimum of $f$. Given $x\in \mathbb{R}^n$ and $S \subset \mathbb{R}^n$, let $d(x,S) := \inf \{\|x - y\| : y \in S\}$. According to the Kurdyka-\L{}ojasiewicz inequality \cite[Theorem 14]{bolte2007clarke} (see also \cite[Theorem 4.1]{attouch2010proximal}), there exist $r,\rho>0$ and a strictly increasing concave continuous semi-algebraic function $\psi:[0,\rho) \rightarrow [0,\infty)$ that is continuously differentiable on $(0,\rho)$ with $\psi(0) = 0$ such that $d(0,\partial f(x)) \geqslant 1/\psi'(|f(x) - f(x^*)|)$ for all $x \in B(x^*,r)$ whenever $0<|f(x)-f(x^*)|< \rho$. After possibly reducing $r$, the inequality holds for all $x \in B(x^*,r)$ such that $f(x) \neq f(x^*)$. 

Let $\epsilon \in (0,r/2)$. By the definition of stability (Definition \ref{def:discrete_lyapunov}), there exist $\delta>0$ and $\bar{\alpha}>0$ such that for all $\alpha \in (0,\bar{\alpha}]$, the subgradient method with constant step size $\alpha$ initialized in $B(x^*,\delta)$ has all its iterates in $B(x^*,\epsilon)$. Since $x^*$ is not a local minimum, we can take an initial iterate $x_0$ in $B(x^*,\delta)$ such that $f(x_0)<f(x^*)$. Let $L>0$ denote a Lipschitz constant of $f$ on $B(x^*,r)$ and let $T := (2\epsilon + \delta)L\psi'(f(x^*) - f(x_0))^2>0$.
Consider the differential inclusion
\begin{equation}
\label{eq:ivp_x_0}
     x'(t) \in - \partial f(x(t)), ~~~ \text{for a.e.}~ 
     t \in (0,T),~~~ x(0) = x_0.
\end{equation}
Since $f$ is locally Lipschitz, the set-valued function $-\partial f$ is upper semicontinuous \cite[2.1.5 Proposition (d)]{clarke1990} with nonempty, compact and convex values \cite[2.1.2 Proposition (a)]{clarke1990}. By Lemma \ref{lemma:euler_method}, there exists $\hat{\alpha} \in (0,\bar{\alpha}]$ such that, for all $\alpha \in (0,\hat{\alpha}]$  and for all sequence $(x_k)_{k\in \mathbb{N}}$ generated by the subgradient method with constant step size $\alpha$ and initialized at $x_0$, there exists a solution $x(\cdot)$ to the differential inclusion \eqref{eq:ivp_x_0} for which $\|\bar{x}(t) - x(t)\| \leqslant \epsilon/2$ for all $t \in [0,T]$, where $\bar{x}:[0,T]\rightarrow\mathbb{R}^n$ is the piecewise linear function defined by $\bar{x}(t) := x_k + (t-\alpha k)/\alpha  (x_{k+1}-x_k)$ for all $t\in [\alpha k , \min\{\alpha(k+1),T\}]$ and $k\in \{0,\hdots,\lfloor T/\alpha\rfloor\}$.

Let us fix some $\alpha \in (0,\min\{\hat{\alpha},T/10\}]$ from now on. Consider a sequence $(x_k)_{k\in\mathbb{N}}$ generated by the subgradient method with constant step size $\alpha$ and initialized at $x_0$. Consider also the linear interpolation $\bar{x}(\cdot)$ of those iterates up to iteration $K+1$ where $K := \lfloor T/\alpha \rfloor$, as well as a solution $x(\cdot)$ to the differential inclusion \eqref{eq:ivp_x_0} such that $\|\bar{x}(t)-x(t)\| \leqslant \epsilon/2$ for all $t\in [0,T]$. Since $f$ is semi-algebraic, by \cite[Lemma 5.2]{davis2020stochastic} (see also \cite{drusvyatskiy2015curves}) it holds that
\begin{equation}
\label{eq:lyapunov}
    f(x(t)) - f(x(0)) = - \int_0^t d(0,\partial f(x(\tau)))^2 d\tau, ~~~ \forall t \in [0,T].
\end{equation}
As a result, $f(x(t)) \leqslant f(x(0)) = f(x_0)$ for all $t \in [0,T]$. Also, $\|x(t) - x^*\| \leqslant \|x(t) - \bar{x}(t)\| + \|\bar{x}(t) - x^*\| \leqslant \epsilon/2 + \epsilon = 3\epsilon/2 \leqslant r$, where the inequality $\|\bar{x}(t) - x^*\| \leqslant \epsilon$ follows from the convexity of $B(x^*,\epsilon)$. Hence $0<f(x^*) - f(x_0) \leqslant f(x^*) - f(x(t))$ and $d(0,\partial f(x(t))) \geqslant  1/\psi'(f(x^*)-f(x(t))) \geqslant 1/\psi'(f(x^*)-f(x_0))$ for all $t\in[0,T]$ by concavity of $\psi$. Together with \eqref{eq:lyapunov}, it follows that $f(x(K\alpha)) - f(x_0) \leqslant -K\alpha M^2$ where $K := \lfloor T/\alpha \rfloor$ and $M:= 1/\psi'(f(x^*)-f(x_0))$. Recall that $L$ is a Lipschitz constant of $f$ on $B(x^*,r)$, so that we have $|f(x(K\alpha))-f(x_0)| \leqslant L \|x(K\alpha)-x_0\|$. We thus obtain the lower bound $\|x(K\alpha)-x_0\| \geqslant K\alpha M^2/L = \lfloor T/\alpha \rfloor \alpha M^2/L \geqslant (T-\alpha)M^2/L \geqslant (T-T/10)M^2/L = 9/10( 2\epsilon + \delta) $ (recall that $T = (2\epsilon + \delta)L/M^2$). Hence $\|x(K\alpha)-x^*\| \geqslant \|x(K\alpha)-x_0\| - \|x_0-x^*\| \geqslant 9/10(2\epsilon + \delta) - \delta = 9 \epsilon/5 - \delta/10 \geqslant 9 \epsilon/5 - \epsilon/10 = 17\epsilon/10$ and finally $\|x_K-x^*\| \geqslant \|x(K\alpha)-x^*\| -\|x_K-x(K\alpha)\| \geqslant 17\epsilon/10 - \epsilon/2 = 6\epsilon/5>\epsilon$. However, since $0<\alpha \leqslant \hat{\alpha} \leqslant \bar{\alpha}$ and $x_0 \in B(x^*,\delta)$, by stability of $x^*$ we have $x_K \in B(x^*,\epsilon)$. We have reached a contradiction. 
\end{proof}

By appealing to Lemma \ref{lemma:euler_method} again, we obtain the following sufficient condition for stability.

\begin{theorem}
\label{thm:sufficient_strict}
Strict local minima of locally Lipschitz semi-algebraic functions are stable.
\end{theorem}
\begin{proof}
Let $x^*$ denote a strict local minimum of locally Lipschitz semi-algebraic function $f:\mathbb{R}^n\rightarrow\mathbb{R}$. Since $x^*$ is a strict local minimum, by the \L{}ojasiewicz inequality \cite[Theorem 0]{kurdyka1998gradients} (see also \cite[\S 2]{lojasiewicz1958}, \cite[\S 17]{lojasiewicz1959}, \cite[(2.1)]{hormander1958division}) and the Kurdyka-\L{}ojasiewicz inequality \cite[Theorem 14]{bolte2007clarke}, there exist $r, \rho>0$ and strictly increasing continuous semi-algebraic functions $\sigma,\psi:[0,\rho)\rightarrow[0,\infty)$ that are continuously differentiable on $(0,\rho)$ with $\sigma(0) = \psi(0) = 0$ such that $\psi$ is concave, $f(x)-f(x^*) \geqslant  \sigma(\|x-x^*\|)$, and $d(0,\partial f(x)) \geqslant 1/\psi'(f(x) - f(x^*))$ for all $x \in B(x^*,r)$ whenever $0<f(x)-f(x^*)< \rho$. After possibly reducing $r$, the two inequalities above hold for all $x \in B(x^*,r)\setminus \{x^*\}$. In order to prove stability, it suffices to prove the statement in Definition \ref{def:discrete_lyapunov} for all $\epsilon>0$ sufficiently small. We may thus restrict ourselves to the case where $0<\epsilon < r$. Given such a fixed $\epsilon$, we next describe a possible choice for $\delta$. 

Given $\Delta \geqslant f(x^*)$, let $L_{x^*}(f,\Delta)$ denote the connected component of the sublevel set $L(f,\Delta) := \{ x \in \mathbb{R}^n : f(x) \leqslant \Delta \}$ containing $x^*$. By taking $\Delta_\epsilon:= f(x^*) + \sigma(\epsilon/2)$, we find that $L_{x^*}(f,\Delta_\epsilon)$ is contained in $B(x^*,\epsilon/2)$. Indeed, one can reason by contradiction and assume that there exists $x \in L_{x^*}(f,\Delta_\epsilon) \setminus B(x^*,\epsilon/2)$. Then $x\notin B(x^*,r)$, otherwise $\sigma(\|x-x^*\|) \leqslant f(x)-f(x^*) \leqslant \sigma(\epsilon/2)$ and thus $\|x-x^*\| \leqslant \epsilon/2$. Therefore $L_{x^*}(f,\Delta_\epsilon)$ is the disjoint union of $L_{x^*}(f,\Delta_\epsilon) \cap \mathring{B}(x^*,r)$ and $L_{x^*}(f,\Delta_\epsilon) \setminus B(x^*,r)$, both of which are nonempty and open in $L_{x^*}(f,\Delta_\epsilon)$. This contradicts the connectedness of $L_{x^*}(f,\Delta_\epsilon)$, which yields that $L_{x^*}(f,\Delta_\epsilon) \subset B(x^*,\epsilon/2)$. By continuity of $f$ we may choose $\delta>0$ such that 

\begin{equation}
\label{eq:inclusions}
    B(x^*,\delta) ~ \subset ~ L_{x^*}(f,\Delta_\epsilon) ~ \subset ~ B(x^*,\epsilon/2).
\end{equation}

We next describe a possible choice for $\bar{\alpha}$. Let $L>\sup\{\|s\|:s\in \partial f(x), x\in B(x^*,\epsilon)\}$ be a Lipschitz constant of $f$ on $B(x^*,\epsilon)$ and let $T := \epsilon/(3L)$, where the supremum is finite due to \cite[Proposition 3 p. 42]{aubin1984differential}. Consider the differential inclusion
\begin{equation}
\label{eq:ivp:levelset}
     x'(t) \in - \partial f(x(t)), ~~~ \text{for a.e.}~ 
     t \in (0,T), ~~~ x(0) \in L_{x^*}(f,\Delta_\epsilon).
\end{equation}
Since $f$ is continuous, the set of initial values $L_{x^*}(f,\Delta_\epsilon)$ is closed. By virtue of the second inclusion in \eqref{eq:inclusions}, $L_{x^*}(f,\Delta_\epsilon)$ is in fact a compact set. Let $\alpha \in (0,T/2]$ and consider a sequence $(x_k)_{k\in\mathbb{N}}$ generated by the subgradient method with constant step size $\alpha$ and initialized in $L_{x^*}(f,\Delta_\epsilon)$. According to the second inclusion in \eqref{eq:inclusions}, the initial iterate $x_0$ lies in $B(x^*,\epsilon/2)$. Hence $\|x_1-x^*\| = \|x_0 - \alpha s_0 - x^*\| \leqslant \|x_0-x^*\| + \|\alpha s_0\| \leqslant \epsilon/2 + L \alpha$ for some $s_0 \in \partial f(x_0)$. Repeating this process until iteration $K+1$ where $K:=\lfloor T/\alpha\rfloor$, we find that $\|x_k-x^*\| \leqslant \epsilon/2 + kL\alpha \leqslant \epsilon/2 + (K+1) L\alpha \leqslant \epsilon/2 + ( T/\alpha+1) L\alpha = \epsilon/2 + (\epsilon/(3L\alpha)+1) L\alpha = 5\epsilon/6+L\alpha \leqslant 5\epsilon/6+LT/2 = 5\epsilon/6+L\epsilon/(3L)/2 = \epsilon$. In other words, the iterates $x_0,\hdots,x_{K+1}$ lie in $B(x^*,\epsilon)$. Let $\epsilon' := \min\{\epsilon L,\sigma(\epsilon/2),\xi^2 T\}/(2L)>0$ where $\xi:= 1/\psi'(\sigma(\epsilon/2)/2)>0$. All of the conditions of Lemma \ref{lemma:euler_method} are met, hence there exists $\bar{\alpha} \in (0,T/2]$ such that, for all $\alpha \in (0,\bar{\alpha}]$ and for all sequence $(x_k)_{k\in\mathbb{N}}$ generated by the subgradient method with constant step size $\alpha$ and initialized in $L_{x^*}(f,\Delta_\epsilon)$, there exists a solution to the differential inclusion \eqref{eq:ivp:levelset} for which $\|\bar{x}(t) - x(t)\| \leqslant \epsilon'$ for all $t\in [0,T]$ where $\bar{x}:[0,T]\rightarrow\mathbb{R}^n$ is the piecewise linear function defined by $\bar{x}(t) := x_k + (t-\alpha k)/\alpha  (x_{k+1}-x_k)$ for all $t\in [\alpha k , \min\{\alpha(k+1),T\}]$ and $k\in \{0,\hdots,\lfloor T/\alpha\rfloor\}$. In particular, it holds that
\begin{equation}
\label{eq:uniform_bound_iterates}
    \|x_k - x(k\alpha)\| \leqslant \epsilon',~~~ k = 0,\hdots,\lfloor T/\alpha \rfloor.
\end{equation}

Having chosen $\delta$ and $\bar{\alpha}$, let us fix some $\alpha \in (0,\bar{\alpha}]$ from now on. Consider a sequence $(x_k)_{k\in\mathbb{N}}$ generated by the subgradient method with constant step size $\alpha$ and initialized in $B(x^*,\delta)$. Our goal is to show that all the iterates lie in $B(x^*,\epsilon)$. According to the first inclusion in \eqref{eq:inclusions}, the initial iterate $x_0$ lies in $L_{x^*}(f,\Delta_\epsilon)$. A previous argument shows that the iterates $x_0,\hdots,x_K$ lie in $B(x^*,\epsilon)$ where $K:=\lfloor T/\alpha \rfloor$. In order to show that the ensuing iterates also lie in $B(x^*,\epsilon)$, we will show that $x_K \in L_{x^*}(f,\Delta_\epsilon)$. The same argument used previously then yields that $x_{K+1},\hdots,x_{2K} \in B(x^*,\epsilon)$. Since $K = \lfloor T/\alpha \rfloor \geqslant \lfloor T/\bar{\alpha}  \rfloor \geqslant 2$, we may conclude by induction that all the iterates belong to $B(x^*,\epsilon)$. This is illustrated in Figure \ref{fig:inclusions}.

\begin{figure}[ht]
\centering
    \begin{tikzpicture}[scale=2.4]
    \draw[blue,opacity=0.2] (4.4,0.4) -- (5.5,.508);
    \fill[gray!50] (5.5,.9) circle (38pt);
    \draw[line width=.2mm] (5.5,.9) circle (38pt);
    \draw[line width=.2mm, name path=A] plot [smooth,tension=0.689] coordinates {(4.4,0.4) (5.1,0.1) (5.8,0.48) (6.4967,0.5682)};
    \draw[line width=.2mm, name path=B] plot [smooth,tension=0.689] coordinates {(6.4967,0.5682)(6.5829,1.0505)(6.0727,1.4115) (5.3159,1.4086) (4.6,1) (4.4,0.4)};
    \tikzfillbetween[of=A and B]{blue, opacity=0.2};
    \fill (4.87,0.79)  node {\normalsize $B(x^*,\delta)$};
    \fill (5.9,1.57)  node {\normalsize $L_{x^*}(f,\Delta_\epsilon)$};
    \draw[line width=.2mm] (5.5,.9) circle (10pt);
    \fill (4.6,2.2)  node {\normalsize $B(x^*,\epsilon/2)$};
    \filldraw (5.5,.9) circle (.2pt);
    \fill (5.61,.92)  node {\normalsize $x^*$};
    \draw[magenta,line width=.3mm] plot [smooth, tension=0.689] coordinates {(6.55,.95) (6.4,.85) (6.2918,.639)};
    \draw[magenta,line width=.3mm] plot [smooth, tension=0.689] coordinates {(6.3,.64) (6.1,.69) (5.85,.7)  (5.6,.6) (5.4,.64)};
    \draw[dashed,line width=.2mm] plot [smooth, tension=0.689] coordinates {(6.6,.8) (6.49,.69) (6.42,.44)};
    \draw[dashed,line width=.2mm] plot [smooth, tension=0.689] coordinates {(6.42,.44) (6.1,.51) (5.85,.52)  (5.6,.42) (5.4,.44)};
    \draw[dashed,line width=.2mm] plot [smooth, tension=0.689] coordinates {(6.5,1.1) (6.32,1) (6.2,.84)};
    \draw[dashed,line width=.2mm] plot [smooth, tension=0.689] coordinates {(6.2,.84) (6.1,.87) (5.85,.88)  (5.6,.78) (5.4,.82)};
    \filldraw[yellow] (6.5,.87) circle (.5pt);
    \filldraw[yellow] (6.33,.84) circle (.5pt);
    \filldraw[yellow] (6.37,.667) circle (.5pt);
    \filldraw[yellow] (6.38,.5) circle (.5pt);
    \filldraw[yellow] (6.23,.52) circle (.5pt);
    \filldraw[yellow] (6.06,.62) circle (.5pt);
    \filldraw[yellow] (5.88,.74) circle (.5pt);
    \filldraw[yellow] (5.65,.68) circle (.5pt);
    \filldraw[yellow] (5.5,.51) circle (.5pt);
    \fill (7.03,.803)  node {\normalsize $x_K$};
    \draw[->] (6.9,.82) -- (6.54,.864);
    \fill (5.0,.385)  node {\normalsize $x_{2K}$};
    \draw[->] (5.15,.405) -- (5.45,.493);
    \end{tikzpicture}
\caption{Induction step}
\label{fig:inclusions}
\end{figure}
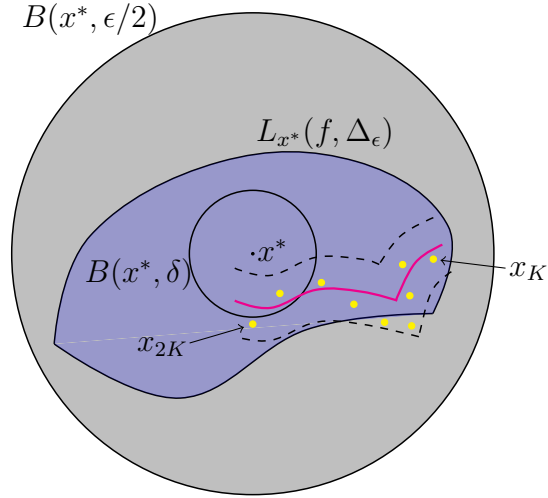

For the remainder of the proof, we seek to show that $x_K \in L_{x^*}(f,\Delta_\epsilon)$. In order to do so, we prove that  $B(x(K\alpha),\epsilon')$ is a connected subset of $L(f,\Delta_\epsilon)$ that has nonempty intersection with $L_{x^*}(f,\Delta_\epsilon)$. Since $L_{x^*}(f,\Delta_\epsilon)$ is a connected component of $L(f,\Delta_\epsilon)$, by maximality and \eqref{eq:uniform_bound_iterates} we then have $x_K \in B(x(K\alpha),\epsilon') \subset L_{x^*}(f,\Delta_\epsilon)$.

We begin by showing that $x(K\alpha) \in B(x(K\alpha),\epsilon') \cap L_{x^*}(f,\Delta_\epsilon)$. Since $f$ is semi-algebraic, by \cite[Lemma 5.2]{davis2020stochastic} (see also \cite{drusvyatskiy2015curves}) it holds that
\begin{equation}
\label{eq:dist}
    f(x(t)) - f(x(0)) = - \int_0^t d(0,\partial f(x(\tau)))^2 d\tau, ~~~ \forall t \in [0,T].
\end{equation}
As a result, $f(x(t)) \leqslant f(x(0)) \leqslant \Delta_\epsilon$ for all $t \in [0,T]$. We thus have that $x(0) \in L_{x^*}(f,\Delta_\epsilon)\cap x([0,T])$, both of which are connected subsets of the sublevel set $L(f,\Delta_\epsilon)$. By maximality of $L_{x^*}(f,\Delta_\epsilon)$, it follows that 
$x(K\alpha)\in x([0,T]) \subset L_{x^*}(f,\Delta_\epsilon)$. 

We next show that $B(x(K\alpha),\epsilon') \subset L(f,\Delta_\epsilon)$. For all $\tilde{x} \in B(x(K\alpha),\epsilon')$, we have
\begin{subequations}
    \begin{align}
        	f(\tilde{x}) - f(x^*) & = f(\tilde{x}) - f(x(K\alpha)) + f(x(K\alpha)) - f(x^*) \label{level_a} \\
	        & \leqslant L\|\tilde{x} - x(K\alpha)\| + \max\{\sigma(\epsilon/2)/2,\sigma(\epsilon/2) - \xi^2 T/2\} \label{level_b} \\
	        & \leqslant L\epsilon' + \sigma(\epsilon/2) - \min\{\sigma(\epsilon/2)/2,\xi^2 T/2\} \label{level_c} \\
            & \leqslant \sigma(\epsilon/2). \label{level_d} 
    \end{align}
\end{subequations}
Indeed, $\tilde{x}$ and $x(K\alpha)$ belong to $B(x^*,\epsilon)$ so we may invoke the Lipschitz constant $L$ of $f$ on $B(x^*,\epsilon)$ in order to bound the first term in \eqref{level_a}. Recall from the previous paragraph that $x(K\alpha)\in L_{x^*}(f,\Delta_\epsilon) \subset B(x^*,\epsilon/2)$ and, since $\epsilon' \leqslant \epsilon/2$, we have $\tilde{x} \in B(x(K\alpha),\epsilon') \subset B(x^*,\epsilon)$. As for the second term in \eqref{level_a}, if it is greater than or equal to $\sigma(\epsilon/2)/2$, then for all $t \in [0,K\alpha]$, we have $\sigma(\epsilon/2)/2 \leqslant  f(x(K\alpha)) - f(x^*) \leqslant  f(x(t)) - f(x^*)$ and thus $d(0,\partial f(x(t))) \geqslant 1/\psi'(f(x(t)) - f(x^*)) \geqslant 1/\psi'(\sigma(\epsilon/2)/2) = \xi$. By \eqref{eq:dist}, it follows that $f(x(K\alpha)) - f(x^*) \leqslant f(x(0)) - f(x^*) - \int_{0}^{K\alpha} \xi^2 d\tau \leqslant f(x(0)) - f(x^*) - K\alpha \xi^2 \leqslant \sigma(\epsilon/2) - \xi^2 T/2$. The last inequality is due to the fact that $x(0) \in L_{x^*}(f,\Delta_\epsilon)$ and $K\alpha = \lfloor T/\alpha\rfloor \alpha \geqslant T-\alpha \geqslant T-\bar{\alpha} \geqslant T/2$. In \eqref{level_c}, we use the fact that $\tilde{x} \in B(x(K\alpha),\epsilon')$ and rewrite the maximum into a minimum. Finally, \eqref{level_d} holds because $\epsilon' \leqslant \min\{\sigma(\epsilon/2),\xi^2 T\}/(2L)$. 
\end{proof}

\begin{figure}[ht]
\centering
\begin{subfigure}{.49\textwidth}
  \centering
  \includegraphics[width=1\textwidth]{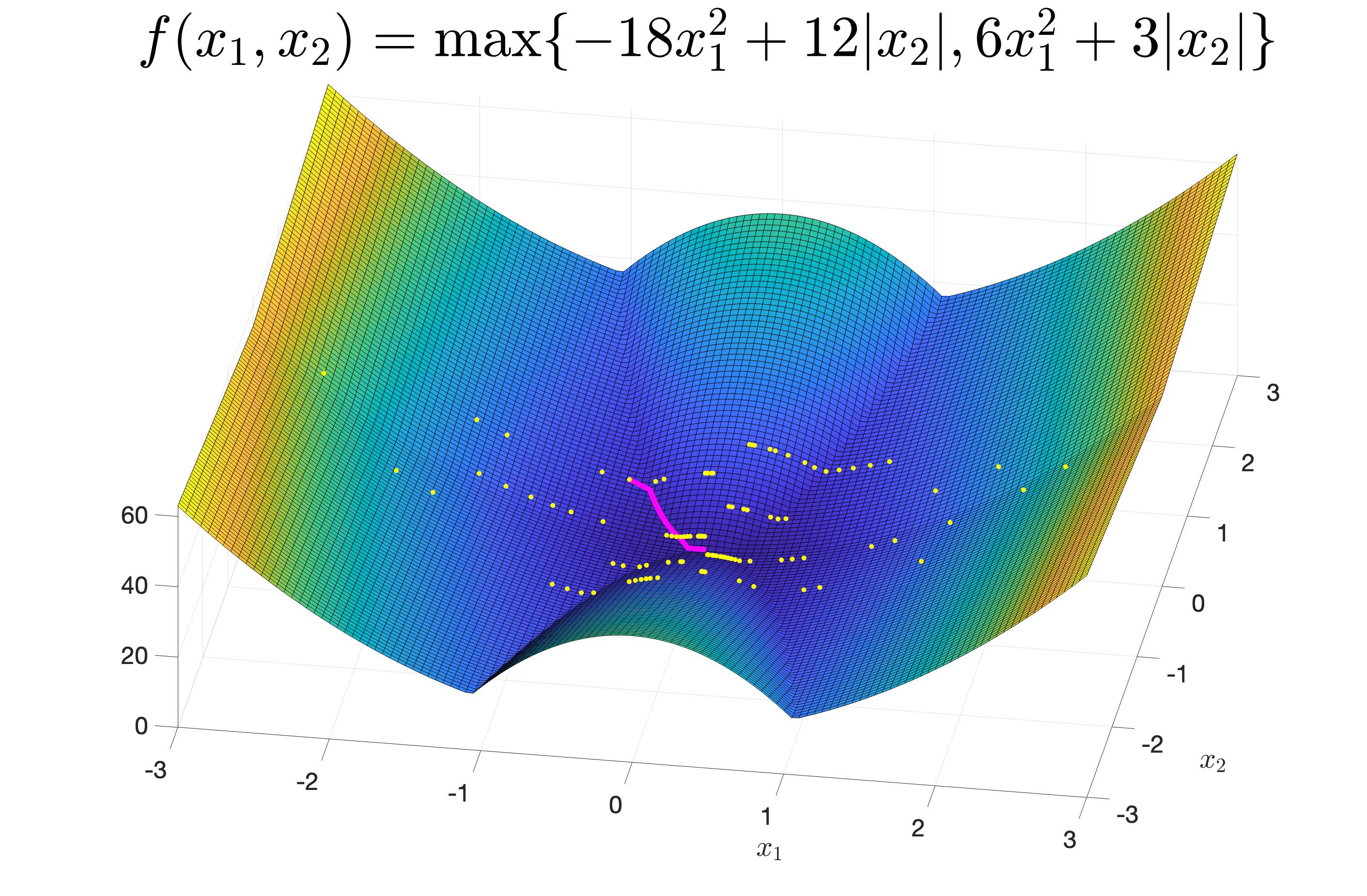}
  \caption{Stable local minimum.}
  \label{fig:stable}
\end{subfigure}
\begin{subfigure}{.49\textwidth}
  \centering
  \includegraphics[width=1\textwidth]{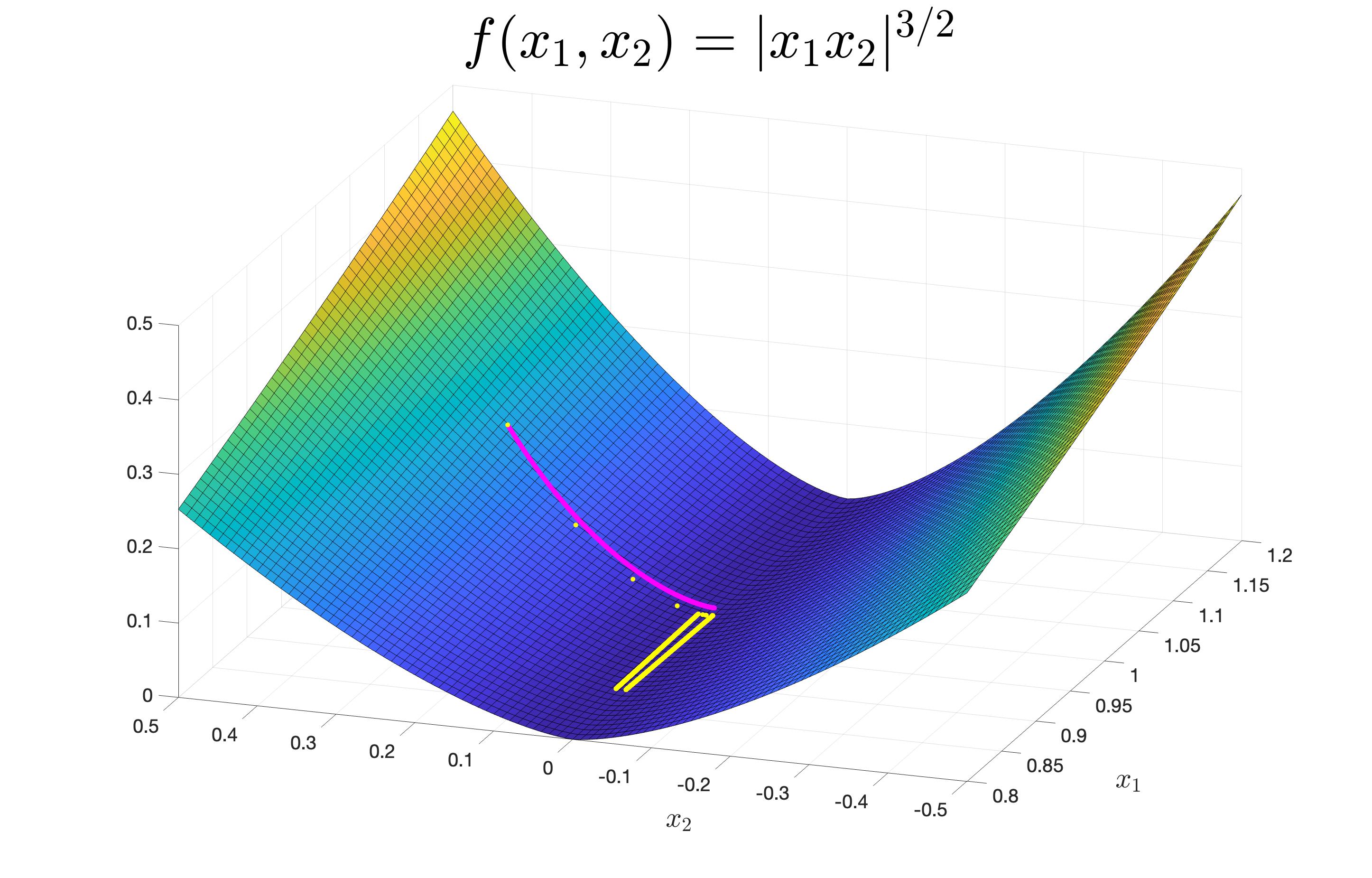}
  \caption{Unstable local minimum.}
  \label{fig:unstable}
\end{subfigure}
\caption{Continuous and discrete subgradient trajectories in magenta and yellow respectively.}
\end{figure}

Observe that non-strict local minima need not be stable. While the strict local minimum in Figure \ref{fig:stable} is stable by Theorem \ref{thm:sufficient_strict}, the non-strict local minimum in Figure \ref{fig:unstable} is unstable. In the former, both continuous and discrete subgradient dynamics are stable, with the continuous trajectory converging to the local minimum, while the discrete trajectory hovers around it. In the latter, continuous and discrete trajectories become decoupled as they approach the local minimum; the continuous dynamics are stable but the discrete dynamics are not.

We conclude this note by proving that the local minimum in Figure \ref{fig:unstable} is unstable with respect to the Euclidean inner product. We show that there exists $\epsilon>0$ such that for all but finitely many constant step sizes $\alpha>0$ and for almost every initial point in $B(x^*,\epsilon)$, at least one of the iterates of the subgradient method does not belong to $B(x^*,\epsilon)$. Let $\epsilon\in (0,1/2]$ and consider the set $S:=\{(x_1,x_2) \in \mathbb{R}^2:x_1x_2=0\}$. By the cell decomposition theorem \cite[(2.11) p. 52]{van1998tame} and \cite[Claim 3]{bolte2020mathematical}, there exist $\alpha_1,\hdots,\alpha_m>0$ such that for all constant step sizes $\alpha \in (0,\infty) \setminus \{\alpha_1,\hdots,\alpha_m\}$, there exists a null subset $I_\alpha \subset \mathbb{R}^2$ such that, for every initial point $(x_1^0,x_2^0) \in \mathbb{R}^2 \setminus I_\alpha$, none of the iterates $(x_1^k,x_2^k)_{k\in \mathbb{N}}$ of the subgradient method belong to the semi-algebraic null set $S$. In this case, the update rule of the subgradient method is given for all $k\in \mathbb{N}$ by
\begin{subequations}
    \begin{align*}
        x_1^{k+1} & = x_1^k - \frac{3}{2} \alpha |x_1^k|^{1/2}|x_2^k|^{3/2} \mathrm{sign}(x_1^k), \\[1mm]
        x_2^{k+1} & = x_2^k - \frac{3}{2} \alpha |x_1^k|^{3/2}|x_2^k|^{1/2} \mathrm{sign}(x_2^k),
    \end{align*}
\end{subequations}
where $\mathrm{sign}(x):=-1$ if $x<0$ and $\mathrm{sign}(x):=1$ if $x>0$. Assume that $(x_1^k,x_2^k) \in B((1,0),\epsilon)$ for all $k\in \mathbb{N}$. Since $\epsilon\leqslant 1/2$, we have $x_1^k \geqslant 1/2$. If $0<|x_2^k|\leqslant \alpha^2/32$ for some $k\in \mathbb{N}$, then $|x_2^{k+1}| = |x_2^k - 3\alpha/2 |x_1^k|^{3/2}|x_2^k|^{1/2} \mathrm{sign}(x_2^k)| \geqslant 3\alpha/2 |x_1^k|^{3/2}|x_2^k|^{1/2} - |x_2^k| \geqslant (3 \alpha /\sqrt{32|x_2^k|} - 1 )|x_2^k| \geqslant 2|x_2^k|$. As a result, $x_2^k$ does not converge to zero. This yields the following contradiction:
\begin{equation*}
    \frac{1}{2} \leqslant x_1^{k+1} = x_1^0 - \frac{3\alpha}{2} \sum_{i=0}^{k} |x_1^i|^{1/2}|x_2^i|^{3/2} \leqslant x_1^0 - \frac{3\alpha}{2\sqrt{2}} \sum_{i=0}^{k} |x_2^i|^{3/2} \rightarrow -\infty.
\end{equation*}

\noindent\textbf{Acknowledgements} We are grateful to the reviewers and editors for their precious time and valuable feedback. We thank Assen Dontchev for fruitful discussions, as well as Salar Fattahi and Richard Zhang for their comments.

\bibliographystyle{abbrv}    
\bibliography{mybib}

\begin{thebibliography}{10}

\bibitem{absil2006stable}
P.-A. Absil and K.~Kurdyka.
\newblock On the stable equilibrium points of gradient systems.
\newblock {\em Systems \& control letters}, 55(7):573--577, 2006.

\bibitem{absil2005convergence}
P.-A. Absil, R.~Mahony, and B.~Andrews.
\newblock Convergence of the iterates of descent methods for analytic cost
  functions.
\newblock {\em SIAM Journal on Optimization}, 16(2):531--547, 2005.

\bibitem{attouch2010proximal}
H.~Attouch, J.~Bolte, P.~Redont, and A.~Soubeyran.
\newblock {Proximal alternating minimization and projection methods for
  nonconvex problems: An approach based on the Kurdyka-{\L}ojasiewicz
  inequality}.
\newblock {\em Mathematics of operations research}, 35:438--457, 2010.

\bibitem{aubin1984differential}
J.-P. Aubin and A.~Cellina.
\newblock {\em Differential inclusions: set-valued maps and viability theory},
  volume 264.
\newblock Springer-Verlag, 1984.

\bibitem{benaim2005stochastic}
M.~Bena{\"\i}m, J.~Hofbauer, and S.~Sorin.
\newblock Stochastic approximations and differential inclusions.
\newblock {\em SIAM Journal on Control and Optimization}, 44(1):328--348, 2005.

\bibitem{bertsekas2015convex}
D.~Bertsekas.
\newblock {\em Convex optimization algorithms}.
\newblock Athena Scientific, 2015.

\bibitem{bianchi2019constant}
P.~Bianchi, W.~Hachem, and A.~Salim.
\newblock Constant step stochastic approximations involving differential
  inclusions: stability, long-run convergence and applications.
\newblock {\em Stochastics}, 91(2):288--320, 2019.

\bibitem{bianchi2022convergence}
P.~Bianchi, W.~Hachem, and S.~Schechtman.
\newblock Convergence of constant step stochastic gradient descent for
  non-smooth non-convex functions.
\newblock {\em Set-Valued and Variational Analysis}, pages 1--31, 2022.

\bibitem{bochnak2013real}
J.~Bochnak, M.~Coste, and M.-F. Roy.
\newblock {\em Real algebraic geometry}, volume~36.
\newblock Springer Science \& Business Media, 2013.

\bibitem{bolte2007clarke}
J.~Bolte, A.~Daniilidis, A.~Lewis, and M.~Shiota.
\newblock Clarke subgradients of stratifiable functions.
\newblock {\em SIAM Journal on Optimization}, 18(2):556--572, 2007.

\bibitem{bolte2020conservative}
J.~Bolte and E.~Pauwels.
\newblock Conservative set valued fields, automatic differentiation, stochastic
  gradient methods and deep learning.
\newblock {\em Mathematical Programming}, pages 1--33, 2020.

\bibitem{bolte2020mathematical}
J.~Bolte and E.~Pauwels.
\newblock A mathematical model for automatic differentiation in machine
  learning.
\newblock {\em NeurIPS}, 2020.

\bibitem{charisopoulos2021low}
V.~Charisopoulos, Y.~Chen, D.~Davis, M.~D{\'\i}az, L.~Ding, and
  D.~Drusvyatskiy.
\newblock Low-rank matrix recovery with composite optimization: good
  conditioning and rapid convergence.
\newblock {\em Foundations of Computational Mathematics}, pages 1--89, 2021.

\bibitem{clarke1990}
F.~H. Clarke.
\newblock {\em Optimization and Nonsmooth Analysis}.
\newblock SIAM Classics in Applied Mathematics, 1990.

\bibitem{davis2020stochastic}
D.~Davis, D.~Drusvyatskiy, S.~Kakade, and J.~D. Lee.
\newblock Stochastic subgradient method converges on tame functions.
\newblock {\em Foundations of computational mathematics}, 20(1):119--154, 2020.

\bibitem{davis2018subgradient}
D.~Davis, D.~Drusvyatskiy, K.~J. MacPhee, and C.~Paquette.
\newblock Subgradient methods for sharp weakly convex functions.
\newblock {\em Journal of Optimization Theory and Applications},
  179(3):962--982, 2018.

\bibitem{drusvyatskiy2015curves}
D.~Drusvyatskiy, A.~D. Ioffe, and A.~S. Lewis.
\newblock Curves of descent.
\newblock {\em SIAM Journal on Control and Optimization}, 53(1):114--138, 2015.

\bibitem{hormander1958division}
L.~H{\"o}rmander.
\newblock On the division of distributions by polynomials.
\newblock {\em Arkiv f{\"o}r matematik}, 3(6):555--568, 1958.

\bibitem{kurdyka1998gradients}
K.~Kurdyka.
\newblock On gradients of functions definable in o-minimal structures.
\newblock In {\em Annales de l'institut Fourier}, volume~48, pages 769--783,
  1998.

\bibitem{lebourg1979generic}
G.~Lebourg.
\newblock Generic differentiability of lipschitzian functions.
\newblock {\em Transactions of the American Mathematical Society},
  256:125--144, 1979.

\bibitem{liapounoff1907probleme}
A.~Liapounoff.
\newblock Probl{\`e}me g{\'e}n{\'e}ral de la stabilit{\'e} du mouvement.
\newblock In {\em Annales de la Facult{\'e} des sciences de Toulouse:
  Math{\'e}matiques}, volume~9, pages 203--474, 1907.

\bibitem{lojasiewicz1958}
S.~\L{}ojasiewicz.
\newblock Division d'une distribution par une fonction analytique de variables
  réelles.
\newblock {\em Comptes rendus hebdomadaires des séances de l'Académie des
  sciences. Paris}, pages 683--686, 1958.

\bibitem{lojasiewicz1959}
S.~\L{}ojasiewicz.
\newblock Sur le problème de la division.
\newblock {\em Studia Mathematica}, page 87–136, 1959.

\bibitem{pham2016genericity}
T.~S. Pham and H.~H. Vui.
\newblock {\em Genericity in polynomial optimization}, volume~3.
\newblock World Scientific, 2016.

\bibitem{rockafellar1981favorable}
R.~T. Rockafellar.
\newblock {Favorable classes of Lipschitz continuous functions in subgradient
  optimization}.
\newblock {\em IIASA Working Paper}, 1981.

\bibitem{sastry2013nonlinear}
S.~Sastry.
\newblock {\em Nonlinear systems: analysis, stability, and control}, volume~10.
\newblock Springer Science \& Business Media, 2013.

\bibitem{van1998tame}
L.~Van~den Dries.
\newblock {\em Tame topology and o-minimal structures}, volume 248.
\newblock Cambridge university press, 1998.

\bibitem{van1996geometric}
L.~Van~den Dries and C.~Miller.
\newblock Geometric categories and o-minimal structures.
\newblock {\em Duke Mathematical Journal}, 84(2):497--540, 1996.

\end{thebibliography}

\end{document}